\documentclass{amsart}

\usepackage{amsmath,amssymb,amscd,amsfonts}

\newtheorem{theorem}{Theorem}
\newcommand{\bt}{\begin{theorem}}
\newcommand{\et}{\end{theorem}}

\newtheorem{lemma}{Lemma}
\newcommand{\bl}{\begin{lemma}}
\newcommand{\el}{\end{lemma}}

\newtheorem{corollary}{Corollary}
\newcommand{\bc}{\begin{corollary}}
\newcommand{\ec}{\end{corollary}}

\DeclareMathOperator{\card}{\text{card}}

\DeclareMathOperator{\qqand}{\qquad\text{and}\qquad}

\newcommand{\N}{\ensuremath{ \mathbf N }}
\newcommand{\PP}{\ensuremath{ \mathbf P}}
\newcommand{\R}{\ensuremath{\mathbf R}}
\newcommand{\Z}{\ensuremath{\mathbf Z}}

\newcommand{\mbx}{\ensuremath{ \mathbf x}}

\newcommand{\beq}{\begin{equation}}
\newcommand{\eeq}{\end{equation}}
\newcommand{\benum}{\begin{enumerate}}
\newcommand{\eenum}{\end{enumerate}}

\newcommand{\mca}{\ensuremath{ \mathcal A}}
\newcommand{\mcb}{\ensuremath{ \mathcal B}}

\title[Multiplicative representations of integers]{Multiplicative representations of integers and Ramsey's theorem}

\author{Melvyn B. Nathanson}
\address{Lehman College (CUNY), Bronx, New York 10468 USA}
\email{melvyn.nathanson@lehman.cuny.edu}

\thanks{Supported in part by a grant from the PSC-CUNY Research Awards Program}

\subjclass[2010]{11B34, 11B75, 11N25, 11B13, 05C55, 05D10.}

\keywords{Multiplicative representations, representation functions, Erd\H os-Tur\' an conjecture, Ramsey theory.}

\date{\today}

\begin{document}

\begin{abstract}
Let $\mcb = (B_1,\ldots, B_h)$ be an $h$-tuple of sets of positive integers.  
Let $g_{\mcb}(n)$ count the number of representations of $n$ in the form $n = b_1\cdots b_h$, 
where $b_i \in B_i$ for all $i \in \{1,\ldots, h\}$.  
It is proved that $\liminf_{n\rightarrow \infty} g_{\mcb}(n) \geq 2$ 
implies $\limsup_{n\rightarrow \infty} g_{\mcb}(n) = \infty$.
\end{abstract}

\maketitle

\section{Does one solution imply many solutions?}

A general problem:  Suppose that an infinite set of equations has the property 
that every equation, or all but finitely many of the equations, or, perhaps, 
infinitely many of the equations have at least one solution.  
Does some equation have at least two solutions?  Three solutions?  Arbitrarily many solutions? 
Infinitely many solutions?  

If the answer is `no,' then  modify the question.  
Suppose that an infinite set of equations has the property 
that every equation, or all but finitely many of the equations, or, perhaps, 
infinitely many of the equations have at least two solutions.  
Does some equation have at least three solutions?  Four solutions?  Arbitrarily many solutions? 
Infinitely many solutions?  

And so on.

This paper considers a problem  of this kind in multiplicative number theory.  

\section{Additive bases}

Let $\N = \{1,2,3,\ldots\}$ be the set of positive integers and $\N_0 = \N \cup \{ 0 \} = \{0,1,2,3,\ldots, \}$ 
the set of nonnegative integers.  
Let $\PP =  \{2,3,5,7,\ldots \}$ be the set of prime numbers.

 
The set  $A$ of nonnegative integers is an \textit{additive  basis of order $h$} 
if every nonnegative integer can be represented as the sum of $h$ elements in the set.
The set $A$ is an \textit{asymptotic additive basis of order $h$} 
if every sufficiently large  integer is the sum of $h$ elements in the set.

For every positive integer $h$, 
the \emph{additive representation function} $r_{A,h}(n)$ counts the number 
of ordered representations of $n$ as a sum of $h$ elements of $A$: 
\[
r_{A,h}(n) = \card\left\{ (a_1,\ldots, a_h) \in A^h: a_1+\cdots + a_h = n \right\}.
\]
Thus,  $A$ is an additive  basis of order $h$ if $r_{A,h}(n) \geq 1$ for all $n \geq 0$, 
and  an asymptotic additive basis of order $h$ if $r_{A,h}(n) \geq 1$ for all $n \geq n_0$.   
In 1941, Erd\H os and Tur\' an~\cite{erdo-tura41} asked if the representation function 
of every asymptotic additive basis of order $2$ is unbounded.
Equivalently, is it true that  
\[
\liminf_{n\rightarrow \infty} r_{A,2}(n) \geq 1 \qquad \text{implies}
 \qquad \limsup_{n\rightarrow \infty} r_{A,2}(n) = \infty?
\]
This problem is still unsolved.  More generally, for $h \geq 2$, is it true that  
\[
\liminf_{n\rightarrow \infty} r_{A,h}(n) \geq 1 \qquad \text{implies}
 \qquad \limsup_{n\rightarrow \infty} r_{A,h}(n) = \infty?
\]

\section{Multiplicative bases and multiplicative systems}

Let $B$ be a set of positive integers.  For every positive integer $h$, 
the \emph{multiplicative representation function} $g_{B,h}(n)$ counts the number 
of ordered representations of $n$ as a product of $h$ elements of $B$: 
\[
g_{B,h}(n) = \card\left\{ (b_1,\ldots, b_h) \in B^h: b_1\cdots  b_h = n \right\}.
\]
The set $B$ is a \emph{multiplicative basis of order $h$} 
if $g_{B,h}(n) \geq 1$ for all $n \geq 1$, 
and an \emph{asymptotic multiplicative basis of order $h$} 
if $g_{B,h}(n) \geq 1$ for all sufficiently large integers $n$.

In 1964, Erd\H os~\cite{erdo64a} proved the following multiplicative analogue 
of the Erd\H os-Tur\' an conjecture.  

\bt[Erd\H os]          \label{MultBasis:theorem:Erdos}
The multiplicative representation function of an asymptotic multiplicative basis of $B$ 
of order $2$ is unbounded.
\et

Equivalently, 
\[
\liminf_{n\rightarrow \infty} g_{B,2}(n) \geq 1 \qquad \text{implies}
 \qquad \limsup_{n\rightarrow \infty} g_{B,2}(n) = \infty.
\]

Erd\H os' proof was graph theoretic.  In 1985, Ne{\v{s}}et{\v{r}}il and R\" odl~\cite{nese-rodl85} 
gave a different proof using Ramsey's theorem.  
In Section~\ref{MultBasis:section:NR}, the Ne{\v{s}}et{\v{r}}il-R\" odl method 
is used to prove Theorem~\ref{MultBasis:theorem:Erdos} 
for asymptotic multiplicative bases of order $h$ for all $h \geq 2$.

There is a natural generalization of a multiplicative basis.  
Let $h \geq 2$ and let $\mcb = (B_1, \ldots, B_h)$ be an $h$-tuple of sets of positive integers.  
For every positive integer $h$, consider 
the \emph{multiplicative representation function}
\[
g_{\mcb}(n) = \card\left\{ (b_1,\ldots, b_h) \in B_1 \times \cdots \times B_h: b_1\cdots  b_h = n \right\}.
\]
The $h$-tuple $\mcb$ is a \textit{multiplicative system of order $h$} 
if $g_{\mcb}(n) \geq 1$ for all $n \geq 1$, and an 
\textit{asymptotic multiplicative system of order $h$} 
if $g_{\mcb}(n) \geq 1$ for all sufficiently large integers $n$.

For example, let $\PP = \{2, 3, 5, \ldots \}$ be the set of prime numbers, and let 
\[
\PP = \PP_1 \cup \cdots \cup \PP_h
\]
be a partition of the prime numbers into $h$ pairwise disjoint nonempty sets.
For all $i = \{1,\ldots, h\}$, let $B_i$ be the set of positive integers all of whose prime factors 
are in $\PP_i$, and let $\mcb = (B_1, \ldots, B_h)$.  
Note that $1 \in B_i$ for all $i$.  
The fundamental theorem of arithmetic is equivalent to the statement that 
$g_{\mcb}(n) = 1$ for all $n = 1,2,3,\ldots$. 
Thus,
\[
\liminf_{n\rightarrow \infty} g_{\mcb}(n) = \limsup_{n\rightarrow \infty} g_{\mcb}(n) =1
\]
and the analogue of Erd\H os' theorem on asymptotic multiplicative bases is false for multiplicative systems. 

This example, however, is misleading.  
For $h \geq 2$, let $M(h)$ consist of all pairs $(s,t)$ such that 
\[
s = \liminf_{n\rightarrow \infty} g_{\mcb}(n) 
\qqand 
 t = \limsup_{n\rightarrow \infty} g_{\mcb}(n) 
\]
for some  asymptotic multiplicative system $\mcb = (B_1, \ldots, B_h)$ of order $h$. 
In Section~\ref{MultBasis:section:MBN} we prove the following result of Nathanson~\cite{nath1987c}:
\[
M(h) =  \{ (1,t):t \in \N \cup \{ \infty\} \} \cup \{ (s,\infty):s=2,\ldots, h \}. 
\]
This implies that  
\[
\liminf_{n\rightarrow \infty} g_{\mcb}(n) \geq 2
\qquad \text{implies}
 \qquad 
 \limsup_{n\rightarrow \infty} g_{\mcb}(n) =\infty.  
 \]

\section{An iterated Ramsey's theorem and multiplicative bases} \label{MultBasis:section:NR}

Let $X$ be a set and let $[X]^k$ be the set of all subsets of $X$ of cardinality $k$.
We have $[X]^0 = \{ \emptyset \}$ and $[X]^1 = \left\{ \left\{ x \right\}:x \in X\right\}$. 
Let $[X]^{< \omega}$ be the set of all finite subsets of $X$.  

The infinite pigeonhole principle states that if an infinite set $X$ is the union 
of a finite number of subsets, then at least one of the subsets is infinite.  
This is equivalent to Ramsey's theorem in the case $k=1$.

\bt[Ramsey's theorem]         \label{MultBasis:theorem:Ramsey-infinite}
Let $X$ be a countably infinite set.  
If $k$ and $r_k$ are nonnegative integers and 
\[
[X]^k = C^{k}_0 \cup C^{k}_1 \cup \cdots  \cup C^{k}_{r_k} 
\]
then there exists an infinite subset $X_k$ of $X$ 
and an integer $\varepsilon_k \in \{0, 1, \ldots, r_k \}$ such that 
\[
[X_k]^k \subseteq C^{k}_{\varepsilon_k}.
\]
\et

A standard reference is \textit{Ramsey Theory} 
by Graham, Rothschild, and Spencer~\cite{grah-roth-spen90}.
We use an iterated form of Ramsey's theorem to prove  
that the representation function of a multiplicative basis is unbounded 
(Theorem~\ref{MultBasis:theorem:Erdos-MBN}).

\bt[Iterated Ramsey's theorem]            \label{MultBasis:theorem:Ramsey-iterated}
Let $X$ be a countably infinite set.   
For  every nonnegative integer $k$, let 
$C^k_0, C^k_1, \ldots, C^k_{r_k}$ be (not necessarily pairwise disjoint) subsets of $[X]^k$ such that 
\beq        \label{MultBasis:partition-xk}
[X]^k =  C^k_0 \cup C^k_1 \cup \cdots  \cup C^k_{r_k}.
\eeq
There exist an infinite decreasing sequence of infinite subsets of $X$ 
\beq                       \label{MultBasis:partition-xk-XX}
X \supseteq  X_0 \supseteq X_1 \supseteq \cdots  \supseteq X_k  \supseteq \cdots  
 \supseteq X_n  \supseteq\cdots 
\eeq
and a sequence $(\varepsilon_k)_{k \in \N_0}$ with 
$\varepsilon_k \in \{0,1,\ldots, r_k \}$ such that 
\beq        \label{MultBasis:partition-sk}
[X_k]^k \subseteq C^{k}_{\varepsilon_k}
\eeq
for all $k \in \N_0$.  
\et

It follows from~\eqref{MultBasis:partition-xk-XX} and~\eqref{MultBasis:partition-sk} 
that $[X_n]^k \subseteq [X_k]^k \subseteq C_{\varepsilon_k}^k$ for all $k \in \{0,1,\ldots, n\}$.   

\begin{proof}
The proof is by induction on $k$.  

Let $X_0 = X$.  We have 
\[
[X_0]^0 = [X]^0 = \{ \emptyset \} = C^0_{\varepsilon_0}
\]
for some $\varepsilon_0 \in \{0,1,\ldots, r_0 \}$.  
This proves the theorem for $k=0$.    
The case $k=1$ follows from the infinite pigeonhole principle.  

Let $k \geq 1$ and let $X_k$ be an infinite subset of $X$ that satisfies~\eqref{MultBasis:partition-sk}.   
For $i \in \{0,1,\ldots, r_{k+1} \}$, let  
\[
D^{k+1}_i = [X_k]^{k+1} \cap C^{k+1}_i.  
\]
Because $[X_k]^{k+1} \subseteq [X]^{k+1}$, we have 
\[
[X_k]^{k+1} = D^{k+1}_0 \cup D^{k+1}_1 \cup \cdots  \cup D^{k+1}_{r_{k+1}}. 
\]
By Ramsey's theorem, there is an  infinite subset $X_{k+1}$ of $X_k$ such that 
\[
[X_{k+1}]^{k+1} \subseteq D^{k+1}_{\varepsilon_{k+1}} \subseteq C^{k+1}_{\varepsilon_{k+1}}
\]
for some $\varepsilon_{k+1} \in \{0,1,\ldots, r_{k+1} \}$.  
This completes the proof.  
\end{proof}

\bt                                    \label{MultBasis:theorem:NR} 
Let $h \geq 2$.  
Let $X$ be a countably infinite set, and let $\mca$ 
be a set of finite subsets of $X$.  
For every finite subset $S$ of $X$, 
let $g_{\mca,h} (S)$ count the number of 
$h$-tuples $(A_{1}, \ldots, A_h )$ in $\mca \times \cdots \times \mca$ such that 
\beq                              \label{MultBasis:theorem:NR-1}
A_{i}\cap A_{j} = \emptyset \qquad \text{for $1 \leq i < j \leq h$}
\eeq
and
\beq                                     \label{MultBasis:theorem:NR-2}
S = \bigcup_{i=1}^h A_{i}.  
\eeq  
Let $W$ be an infinite set of positive integers such that $g_{\mca,h} (S) \geq 1$ 
for all $n \in W$ and for all but finitely many sets $S \in [X]^n$.  
For every integer $n$, there is a finite subset $S$ of $X$ with 
$g_{\mca,h}(S) \geq n$. 
\et

\begin{proof}
For all $k \in \N_0$, let 
\[
C_1^k = [X]^k \cap \mca \qqand C_0^k = [X]^k \setminus C_1^k.
\]
Thus, 
\[
[X]^k = C_0^k \cup C_1^k \qqand   C_0^k \cap C_1^k = \emptyset.
\]
We apply Theorem~\ref{MultBasis:theorem:Ramsey-iterated} with $r_k = 1$ for all $k \in \N_0$.  
For every positive integer $n$, there is an infinite subset $X_n$ of $X$ 
such that, for all $k \in \{0,1,\ldots, n\}$, 
 there exists $\varepsilon_k \in \{0,1\}$ with 
 \[
 [X_n]^k \subseteq C_{\varepsilon_k}^k. 
\]
 
Let $n \in W$.  The infinite set $X_n$ contains infinitely many subsets of cardinality $n$.  
Choose a subset $S$ of $X_n$ of cardinality $n$ for which 
there exists an $h$-tuple $(A_{1}, \ldots, A_{h})$ in $\mca \times \cdots \times \mca$ satisfying 
conditions~\eqref{MultBasis:theorem:NR-1} and~\eqref{MultBasis:theorem:NR-2}.  
We have $k_i = |A_i| \in \{0,1,\ldots, n\}$ for $i \in \{1,\ldots, h\}$.  
The sets $A_i$ are pairwise disjoint, and so  
\[
\sum_{i=1}^h k_i = \sum_{i=1}^h |A_i| = |S| = n. 
\]
Because $A_i  \subseteq S \subseteq X_n \subseteq X$ and $A_i \in \mca$, 
we have 
\[
A_i \in   [X_n]^{k_i} \cap \mca \subseteq  [X]^{k_i} \cap \mca  = C_1^{k_i}. 
\]
It follows that $\varepsilon_{k_i} = 1$ and 
\beq                                     \label{MultBasis:theorem:NR-XC}
[X_n]^{k_i} \subseteq C_1^{k_i} \subseteq \mca
\eeq
for all $i \in \{1,\ldots, h \}$.  
Thus, every subset of $X_n$ of cardinality $k_i$ belongs to \mca.  
The multinomial coefficient 
\[
\binom{n}{k_1,\ldots, k_h} = \frac{n!}{k_1! \cdots k_h!} 
\]
counts the number of partitions of $S$ into pairwise disjoint subsets $A'_1,\ldots, A'_h$ 
of cardinalities $k_1,\ldots, k_h$, respectively.  
It follows from~\eqref{MultBasis:theorem:NR-XC} that $A'_i \in \mca$ for all $i \in \{1,\ldots, h\}$, 
and so 
\[
g_{\mca}(S)  \geq \binom{n}{k_1,\ldots, k_h}. 
\] 
If $k_i < n$ for all $i$, then 
\[
g_{\mca}(S) \geq \binom{n}{k_1,\ldots, k_h} \geq n.
\]
If $k_j = n$ for some $j$, then $k_i = 0$ for all $i \neq j$ and 
$\binom{n}{k_1,\ldots, k_h} = 1$. 
In this case, $S \in \mca$ and $[X_n]^n \subseteq C_1^n$.  
Thus, every $n$-element subset of $X_n$ belongs to \mca.  
Also, $\emptyset \in \mca$.  
Let $S'$ be a subset of $X_n$ of cardinality $2n$.  
There are $\binom{2n}{n}$ partitions of $S$ of the form $T = A_1 \cup A_2 \cup \cdots \cup A_h$, 
where $|A_1| = |A_2| = n$ and $|A_i| = 0$ for $i \in \{3,\ldots, h\}$.  
It follows that 
\[
g_{\mca}(S') \geq \binom{2n}{n} \geq n.
\]
This completes the proof.  
\end{proof}

\bt        \label{MultBasis:theorem:Erdos-MBN}
The multiplicative representation function of an asymptotic multiplicative basis $B$ 
of order $h$ is unbounded.
\et

Equivalently, 
\[
\liminf_{n\rightarrow \infty} g_{B,h}(n) \geq 1 \qquad \text{implies}
 \qquad \limsup_{n\rightarrow \infty} g_{B,h}(n) = \infty.
\]

\begin{proof} 
A \textit{square-free (or quadratfrei) integer} is a positive integer that is the product of distinct primes.  
Equivalently, a positive integer $n$ is square-free if it is not divisible by the square of a prime.  
 Let $Q = \{1,2,3,5,6,7,10,11,\ldots\}$ be the set of square-free integers.  
 Define the function $\Phi: Q \rightarrow [\PP]^{< \omega}$ as follows:
 For $q \in Q$, 
 \[
 \Phi(q) = \{p \in \PP: \text{$p$ divides $q$}\}.
 \]
The function $ \Phi$ is a bijection, with inverse function 
$ \Phi^{-1}: [\PP]^{< \omega} \rightarrow Q$ defined as follows: 
For $S \in  [\PP]^{< \omega}$, 
\[
\Phi^{-1}(S) = \prod_{p\in S} p.
\]

If $n$ is square-free and if $b_1,\ldots, b_h$ are positive integers such that $n = b_1\cdots b_h$, 
then the integers $b_1,\ldots, b_h$ are square-free and pairwise relatively prime.  
Let  $S = \Phi(n)$ and let $S_i = \Phi(b_i)$ for $i \in \{1,\ldots, h\}$. 
It follows that $S= \bigcup_{i=1}^h S_i$ and that $S_1,\ldots, S_h$ are pairwise disjoint 
sets in $ [\PP]^{< \omega}$. 

Conversely, let $S_1,\ldots, S_h$ be pairwise disjoint sets in $ [\PP]^{< \omega}$  
and $S = \bigcup_{i=1}^h S_i$.  
Let  $n = \Phi^{-1}(S)$ and let $b_i = \Phi^{-1}(S_i)$ for $i \in \{1,\ldots, h\}$. 
We have 
\[
n = \Phi^{-1}(S) = \prod_{i=1}^h \Phi^{-1}(S_i) = \prod_{i=1}^h b_i. 
\]
Thus, there is a one to one correspondence between partitions of a finite set 
of primes into $h$ pairwise disjoint subsets and representations of a square-free 
integer as a product of $h$ factors.   

Let $B$ be an asymptotic multiplicative basis of order $h$.  Let 
\[
\mca = \{ \Phi^{-1} (b): b \in B \cap Q\} \subseteq [\PP]^{< \omega}.
\]
Because $B$ is an asymptotic multiplicative basis, 
every square-free integer, with a finite number of exceptions,  
has a representation as a product of $h$ pairwise relatively prime square-free integers in $B$.  
Thus, with a finite number of exceptions, every set $S \in [\PP]^{< \omega}$ 
has a representation as the union of $h$ pairwise disjoint sets in \mca.  
Theorem~\ref{MultBasis:theorem:NR} implies that for every integer $n$ there is a set 
$S \in [\PP]^{< \omega}$ with at least $n$ representations as the union 
of $h$ pairwise disjoint sets in \mca.  It follows that the integer $q = \Phi(S) \in Q$ 
has at least $n$ representations as a product of $h$ pairwise relatively prime 
square-free integers in $B$, and so $g_{B,h}(q) \geq n$.  
This completes the proof.  
\end{proof}

\section{A doubly iterated Ramsey's theorem}        \label{MultBasis:section:DoubleRamsey}

We use an iterated form of the iterated Ramsey's theorem to study  
the representation function of a multiplicative system (Theorem~\ref{MultBasis:theorem:MBN}).  
We need the following elementary Boolean identity.

\bl                                 \label{MultBasis:lemma:Boolean}
Let $Y$ and $I$ be sets.  For all $i \in I$, let $J_i$ be a set, and let $\{C_{i,j_i}: j_i \in J_i\}$ 
be a set of subsets of $Y$ such that 
\[
Y = \bigcup_{j_i \in J_i} C_{i,j_i}.
\]
Then
\[
Y = \bigcup_{(j_i)\in \prod_{i\in I} J_i}  \hspace{0.2cm}  \bigcap_{i\in I} C_{i,j_i} 
\]
where $\prod_{i\in I} J_i$ is the Cartesian product of the sets $J_i$.  
\el

\begin{proof}
Because the sets $C_{i,j_i}$ are subsets of $Y$, we have 
\[
\bigcup_{(j_i)\in \prod_{i\in I} J_i} \hspace{0.2cm}  \bigcap_{i\in I} C_{i,j_i} \subseteq  Y.  
\]

Let $y\in Y$.  For all  $i \in I$, there exists $j_i \in J_i$ such that $y \in C_{i,j_i}$ and so 
\[
y \in \bigcap_{i\in I} C_{i,j_i}.
\]
Therefore, 
\[
y \in \bigcup_{(j_i)\in \prod_{i\in I} J_i}  \hspace{0.2cm} \bigcap_{i\in I} C_{i,j_i}
\]
and
\[
Y \subseteq \bigcup_{(j_i)\in \prod_{i\in I} J_i}  \hspace{0.2cm} \bigcap_{i\in I} C_{i,j_i}. 
\]
This completes the proof.  
\end{proof}

\bt[Doubly iterated Ramsey's theorem]            \label{MultBasis:theorem:Ramsey-iterated-iterated}
Let $X$ be a countably infinite set.   
Let $I_k$ be a finite set for every nonnegative integer $k$.  
For all $i \in I_k$, let $J^k_i$ be a finite set and let $\{C^k_{i,j_i} : j_i \in J^k_i\} $ 
be a set of (not necessarily pairwise disjoint) subsets of $[X]^k$ such that  
\beq        \label{MultBasis:partition-xk-111}
[X]^k = \bigcup_{j_i \in J_i} C^k_{i,j_i}. 
\eeq 
There is an infinite decreasing sequence of infinite subsets of $X$ 
\beq        \label{MultBasis:partition-xk-X}
X \supseteq X_0 \supseteq X_1 \supseteq \cdots  \supseteq X_k  \supseteq X_{k+1}  \supseteq\cdots 
\eeq
and, for all $k \in \N_0$ and for all $i \in I_k$, 
there is an integer $\varepsilon_{i,k} \in J_i^k$ such that 
\beq        \label{MultBasis:partition-sk-111}
[X_k]^k \subseteq C^{k}_{i,\varepsilon_{i,k}}. 
\eeq 
\et

It follows from~\eqref{MultBasis:partition-xk-X} and~\eqref{MultBasis:partition-sk-111} 
that $[X_n]^k \subseteq [X_k]^k \subseteq 
C^k_{i, \varepsilon_{i,k}}$  for all $k \in \{0,1,\ldots, n\}$.

\begin{proof}
It suffices to prove that there exists $(\varepsilon_{i,k}) = (\varepsilon_{i,k})_{i\in I_k}  \in \prod_{i\in I_k} J_i^k$ such that 
\[
[X_k]^k \subseteq \bigcap_{i \in I_k } C^{k}_{i,\varepsilon_{i,k}}.
\]
Applying  Lemma~\ref{MultBasis:lemma:Boolean} with $Y = [X]^k$, we obtain 
\[
[X]^k = \bigcup_{(\varepsilon_{i,k})\in \prod_{i\in I_k} J^k_i}  \hspace{0.2cm}  \bigcap_{i\in I} C_{i,\varepsilon_{i,k}} 
\]
The set $\prod_{i\in I_k} J^k_i$ is finite because the sets $I_k$ and $J^k_i$ are finite.  
The doubly iterated Ramsey's theorem follows immediately from the iterated 
Ramsey's theorem (Theorem~\ref{MultBasis:theorem:Ramsey-iterated}).
This completes the proof.  
\end{proof}

\bt                      \label{MultBasis:theorem:key}
Let $h \geq 2$.  
Let $X$ be a countably infinite set, and let $\mca^* = (\mca_1, \ldots, \mca_h)$ 
be an $h$-tuple of sets of finite subsets of $X$.  
For every  subset $S$ of $X$, 
let $g_{\mca^*} (S)$ count the number of 
$h$-tuples $(A_{1}, \ldots, A_h )$ in $\mca_1 \times \cdots \times \mca_h$ such that 
\beq                              \label{MultBasis:SA1}
A_{i}\cap A_{j} = \emptyset \qquad \text{for $1 \leq i < j \leq h$}
\eeq
and                                           
\beq                              \label{MultBasis:SA2}        
S = \bigcup_{i=1}^h A_{i}.  
\eeq  
Let $W$ be an infinite set of positive integers such that $g_{\mca^*} (S) \geq 2$ 
for all $n \in W$ and for all but finitely many sets $S \in [X]^n$.  
For every integer $n$, there is a finite subset $S$ of $X$ with 
$g_{\mca^*}(S) \geq n$.  
\et

\begin{proof}
For all $k \in \N_0$, let  $I_k = \{1,\ldots, h\}$ and, for all $i  \in I_k$, let $J_i^k = \{0,1\}$.  
The sets 
\[
C_{i,1}^k =  [X]^k \cap \mca_i 
 \]
and 
\[
C_{i,0}^k =  [X]^k \setminus C_{1,i}^k 
 \]
partition $[X]^k$.  Thus, 
 \[
[X]^k  = C_{i,0}^k \cup C_{i,1}^k \qqand C_{i,0}^k \cap C_{i,1}^k = \emptyset.
 \]
By Theorem~\ref{MultBasis:theorem:Ramsey-iterated-iterated}, 
there is an infinite decreasing sequence 
\[
X \supseteq X_0 \supseteq X_1 \supseteq \cdots  \supseteq X_k  \supseteq\cdots  \supseteq X_n  \supseteq\cdots 
\]
of infinite subsets of $X$ such that, for all $n \in \N_0$, 
the set $X_n$ has the following property:  
For all $k \in \{ 0,1,\ldots, n\}$ and $i \in I_k$,  
there exists $\varepsilon_{i,k} \in \{0,1\}$ with 
\[
[X_n]^k \subseteq [X_k]^k \subseteq C_{i,\varepsilon_{i,k}}^k.
 \]
If there is a subset $A_i$ of $X_n$ of cardinality $k_i \leq n$ with $A_i \in \mca_i$, then 
$A_i \in C_{i,1}^{k_i}$ and so $[X_n]^{k_i} \subseteq C_{i,1}^{k_i}$, 
that is,  every subset of $X_n$ of cardinality ${k_i}$ is in $\mca_i$.
 
Let $n \in W$.  
Because $X_n$ is an infinite subset of $X$, there is a set $S \in [X_n]^n$ 
with $g_{\mca^*}(S) \geq 2$.  
Thus, there are  distinct  $h$-tuples $(A_1, \ldots, A_h)$ and $(A'_1, \ldots, A'_h)$ 
in $\mca_1 \times \cdots \times \mca_h$
that satisfy conditions~\eqref{MultBasis:SA1} and~\eqref{MultBasis:SA2}.  
Let 
\[
|A_i| = k_i \qqand |A'_i| = k'_i 
\]
for $i \in \{1,\ldots, h\}$.  The sets $A_1,\dots, A_h$ are pairwise disjoint 
and the sets $A'_1,\dots, A'_h$ are pairwise disjoint, and so 
\[
\sum_{i=1}^h k_i = \sum_{i=1}^h k'_i = n.
\] 
Because $A_i \in [X_n]^{k_i} \cap \mca_i$ and $A'_i \in [X_n]^{k'_i} \cap \mca_i$, 
every subset of $X_n$ of cardinality $k_i$ or $k'_i$ is in $\mca_i$  
It follows that if $S = \bigcup_{i=1}^h B_i$ is a partition of $S$ with $|B_i| = k_i$ for 
all $i \in \{1,\ldots, h\}$,  then $B_i \in \mca_i$ for all $i \in \{1,\ldots, h\}$.  
The number of such partitions is the multinomial coefficient  
\[
\binom{n}{k_1,\ldots, k_h}.    
\]
If $k_i < n$ for all $i \in \{1,\ldots, h\}$, then  
\[
\binom{n}{k_1,\ldots, k_h} \geq n.    
\]
Similarly,  if $S = \bigcup_{i=1}^h B'_i$ is a partition of $S$ with $|B'_i| = k'_i$ for 
all $i \in \{1,\ldots, h\}$,  then $B'_i \in \mca_i$ for all $i \in \{1,\ldots, h\}$.  
The number of such partitions is the multinomial coefficient  
\[
\binom{n}{k'_1,\ldots, k'_h}.    
\]
If $k'_i < n$ for all $i \in \{1,\ldots, h\}$, then  
\[
\binom{n}{k'_1,\ldots, k'_h} \geq n.    
\]

If $k_i = k'_j= n$ for some $i,j \in \{1,\ldots, h \}$,   
then $B_i = B'_j= S$ and $B_{\ell} = \emptyset$ for all $\ell \neq i$ and 
$B'_\ell = \emptyset$ for all $\ell \neq j$.  
Because $(B_1,\ldots, B_h) \neq (B'_1,\ldots, B'_h)$, we have $i \neq j$.
Every subset of $X_n$ of cardinality $n$ is in  both $\mca_i$ and $\mca_j$.  
Also, $\emptyset \in \mca_i$ for all $i \in \{1,\ldots, h\}$.  

Let $S_1$ and $S_2$ be disjoint subsets of $[X_n]^n$, and let $S = S_1 \cup S_2$.
We have $|S| = 2n$.  There are $\binom{2n}{n}$ subsets $U$ of $S$ of cardinality $n$.
For each of these sets $U$, let $B_i = U$, $B_j= S\setminus U$, and $B_{\ell} = \emptyset$ 
for $\ell \neq i,j$.  
This gives $\binom{2n}{n} > n$ ordered $h$-tuples 
$(B_1,\ldots, B_h) \in \mca_1 \times \cdots \times \mca_h$ 
that satisfy conditions~\eqref{MultBasis:SA1} and~\eqref{MultBasis:SA2}.  
This completes the proof.  
\end{proof}

\section{Representation functions of multiplicative systems}   \label{MultBasis:section:MBN}

\bt                       \label{MultBasis:theorem:MBN}
For $h \geq 2$, let $M(h)$ consist of all pairs $(s,t)$ such that 
\[
s = \liminf_{n\rightarrow \infty} g_{\mcb}(n) 
\qqand 
 t = \limsup_{n\rightarrow \infty} g_{\mcb}(n) 
\]
for some  asymptotic multiplicative system $\mcb$ of order $h$. 
Then 
\[
M(h) =  \{ (1,t):t \in \N  \cup \{\infty\} \} \cup \{ (s,\infty):s \in \{2,\ldots, h\} \}.
\]
\et

\begin{proof}
For every asymptotic multiplicative system $\mcb = (B_1, \ldots, B_h)$ of order $h$, 
we have 
\[
s = \liminf_{n\rightarrow \infty} g_{\mcb}(n) \geq 1.
\]
Let $p$ be a prime number.  For all $(b_1,\ldots, b_h) \in B_1 \times \cdots \times B_h$, 
we have  $p = b_1\cdots b_h$ if and only if $b_j = p$ for some $j \in\{1,\ldots, h\}$ and $b_i = 1$ for all $i \neq j$.
It follows that $g_{\mcb}(p) \leq h$, and so 
\[
s = \liminf_{n\rightarrow \infty} g_{\mcb}(n)   \leq h. 
\]
For every positive integer $t$, define the sets 
\[
B_1 = \N, \quad B_2 = \{2^k:0 \leq k \leq t-1 \},\quad B_i = \{ 1\}\quad \text{for $i=3,\ldots, h$}
\]
and let  $\mcb = (B_1, \ldots, B_h)$.  
If $n$ is an odd positive integer, then 
\[
g_{\mcb}(n) = 1  
\]
and so 
\[
 \liminf_{n\rightarrow \infty} g_{\mcb}(n) = 1.
\]
Let $n = 2^{\ell -1} m$, where $m$ is odd.  
We have $(b_1, b_2,\ldots, b_h) \in B_1 \times B_2 \times  \cdots \times B_h$ and  $n = b_1 b_2 \cdots b_h$ 
if and only if 
\[
(b_1, b_2, \ldots, b_h) = \left(2^{\ell-1-i} m, 2^i,1,\ldots, 1 \right) 
\]
for some $i \in \{0,1,\ldots, \min(\ell -1,t-1) \} $.
It follows that $g_{\mcb}(n) = \min(\ell,t)$ 
and so 
\[
 \limsup_{n\rightarrow \infty} g_{\mcb}(n) = t.
\]  
Therefore, $(1,t) \in M(h)$ for all $t \in \N$.  

Let $\mcb = (B_1, \ldots, B_h)$, where 
\[
B_1 = \N, \quad B_2 = \{2^k: k \in \N_0\},\quad B_i = \{ 1\}\quad \text{for $i=3,\ldots, h$.}
\]
We have 
\[
  \liminf_{n\rightarrow \infty} g_{\mcb}(n) = 1 \qqand \limsup_{n\rightarrow \infty} g_{\mcb}(n) = \infty 
\]  
and so $(1, \infty) \in M(h)$.

For $s \in \{2,3,\ldots, h\}$, let $\mcb = (B_1, \ldots, B_h)$ be the multiplicative system 
defined by 
\[
B_1 = \N, \qquad B_i = \PP \cup \{1\} \quad \text{for $i=2, \ldots, s$,}
\]
and 
\[
 B_i =  \{1\} \quad \text{for $i=s+1, \ldots, h$.} 
\]
We have 
\[
g_{\mcb}(p) = s \qquad \text{for all $p \in \PP$}
\]
and
\[
g_{\mcb}(n) \geq s \qquad \text{for all $n \geq 2$.}
\]
Therefore,
\[
 \liminf_{n\rightarrow \infty} g_{\mcb}(n) = s.
\]
If $n \in \N$ has $k$ distinct prime factors, then $g_{\mcb}(n) \geq k$ and so 
\[
 \limsup_{n\rightarrow \infty} g_{\mcb}(n) = \infty.
\]
It follows that  
\[
 \{ (s,\infty):s \in \{2,\ldots, h\} \}\subseteq M(h).
\]

We must prove that if $(s,t) \in M(h)$ and $s \geq 2$, then $t = \infty$.  
As in the proof of Theorem~\ref{MultBasis:theorem:Erdos-MBN}, 
we consider the set $Q$ of square-free integers, and the functions 
\[
\Phi: Q \rightarrow [\PP]^{<\omega} \qqand 
\Phi^{-1}:  [\PP]^{<\omega} \rightarrow Q
\]
defined by 
\[
\Phi(q) = \{ p\in \PP:  \text{$p$ divides $q$} \} \qquad \text{for all $q \in Q$} 
\]
and
\[
\Phi^{-1}(S) = \prod_{p\in S} p  \qquad \qquad \text{for all $S \in  [\PP]^{<\omega}$.} 
\]
There is a one to one correspondence between  ordered partitions of a finite set 
of primes into $h$ pairwise disjoint subsets and ordered representations of a square-free 
integer as a product of $h$ factors.    
Let $S_1,\ldots, S_h$ be pairwise disjoint sets in $ [\PP]^{< \omega}$, 
and let $S = \bigcup_{i=1}^h S_i$.
Let $n = \Phi^{-1}(S)$ and $b_i = \Phi^{-1}(S_i)$ for all $i \in \{1,\ldots, h\}$.  
We have 
\[
n = \Phi^{-1}(S) = \prod_{i=1}^h \Phi^{-1}(S_i) = \prod_{i=1}^h b_i. 
\]
Conversely, if $b_1,\ldots, b_h$ are pairwise relatively prime square-free integers 
and $n = b_1\cdots b_h$ with $b_i \in B_i$, 
then the sets $\Phi(b_1),\ldots, \Phi(b_h)$ are pairwise disjoint and 
$\Phi(n) = \bigcup_{i=1}^h \Phi(b_i)$.

For all $i \in \{1,\ldots, h\}$, let  
\[
\mca_i = \{ \Phi(b_i) : b_i \in Q \cap B_i \} \subseteq [\PP]^{<\omega}
\]
where $Q \cap B_i $ is the set of square-free integers in the set $B_i$.
Let $\mca^* = (\mca_1,\ldots, \mca_h)$.  
We have $b_i \in Q \cap B_i$ if and only if $\Phi(b_i) \in \mca_i$, and so 
\[
g_{\mcb}(q) = g_{\mca^*} \left(\Phi(q)\right)
\]
for all square-free integers $q$.  
The inequality $\liminf_{q \in Q} g_{\mcb}(q) \geq 2$ is equivalent to 
$\liminf_{S \in [\PP]^{<\omega}} g_{\mca^*}(S) \geq 2$.  
Theorem~\ref{MultBasis:theorem:key} implies 
$  \limsup_{S \in [\PP]^{<\omega}} g_{\mca^*}(S) = \infty$, and so 
\begin{align*}
 \limsup_{n\rightarrow \infty} g_{\mcb}(n) 
 & = \limsup_{q \in Q} g_{\mcb}(q) \\
 & = \limsup_{q \in Q}  g_{\mca^*}(\Phi(q) )  \\
& =  \limsup_{S \in [\PP]^{<\omega}} g_{\mca^*}(S) \\
& = \infty.  
\end{align*} 
This completes the proof. 
\end{proof}

\bc
For every  asymptotic multiplicative system \mcb\ of order $h$, 
 \[
\liminf_{n\rightarrow \infty} g_{\mcb}(n) \geq 2
\qquad \text{implies} \qquad 
 \limsup_{n\rightarrow \infty} g_{\mcb}(n) = \infty.  
\]
\ec

\section{Open problems}
Theorems~\ref{MultBasis:theorem:Erdos-MBN} and~\ref{MultBasis:theorem:MBN} 
are answers to questions of the following kind.  
Suppose that a representation function has at least one, or two, or three, or \ldots solutions.   
Does this imply that there are occasionally even more  solutions?  
Does this imply that the representation function can be arbitrarily large?  
Here are related open problems.  

\benum
\item
Let \mca\ be a set of nonnegative integers that is an asymptotic additive basis of order 2.
Suppose that every sufficently large integer has at least 1000 representations 
as the sum of two elements of \mca?  
Does some integer have 1010 representations? 

\item
Let $\mca$ be a set of vectors in $\Z^n$.  The set \mca\ is an \textit{inner product basis} 
for a set $W$ of integers if every integer in $W$ is the inner product of two vectors in \mca.  
For what infinite sets $W$ with  inner product basis \mca\  is the representation function bounded?
For what infinite sets $W$ with  inner product basis \mca\  is the representation function unbounded?
 
 \item
Associated to every $n\times n$ matrix $M$ is the quadratic form $Q_M$, 
defined by $Q_M (\mbx) = \mbx^t M \mbx$ for $\mbx \in \R^n$.  
Let $M$ be an $n\times n$ matrix  with nonnegative integral coordinates.  
Let \mca\ be a set vectors with nonnegative  integral coordinates, 
and let $W$ be an infinite 
set of nonnegative integers such that $W \subseteq Q_M(\mca)$.
When is the representation function bounded?
When is the representation function unbounded?

\eenum

\def\cprime{$'$} \def\cprime{$'$} \def\cprime{$'$}
\providecommand{\bysame}{\leavevmode\hbox to3em{\hrulefill}\thinspace}
\providecommand{\MR}{\relax\ifhmode\unskip\space\fi MR }
\providecommand{\MRhref}[2]{%
  \href{http://www.ams.org/mathscinet-getitem?mr=#1}{#2}
}
\providecommand{\href}[2]{#2}

\end{document}